\documentclass[11pt]{article}
\usepackage{amsmath, amsthm}
\usepackage{graphicx}
\usepackage{amssymb}
\usepackage{amsfonts}

\newtheorem{theorem}{Theorem}
\newtheorem{corollary}{Corollary}
\newtheorem{lemma}{Lemma}
\newtheorem{proposition}{Proposition}
\newtheorem*{theorem*}{Theorem}

\theoremstyle{remark}
\newtheorem{remark}{Remark}
\begin{document}

\title{\textbf{Bahadur--Kiefer Representations for Time Dependent Quantile
Processes }}
\author{P\'eter Kevei\thanks{This research was funded by a postdoctoral fellowship
of the Alexander von Humboldt Foundation.} \\
Center for Mathematical Sciences, Technische Universit\"at M\"unchen \\
Boltzmannstra{\ss}e 3, 85748 Garching, Germany, and \\
MTA-SZTE Analysis and Stochastics Research Group \\
Bolyai Institute, Aradi v\'{e}rtan\'{u}k tere 1, 6720 Szeged, Hungary \\
e-mail: \texttt{peter.kevei@tum.de} \smallskip\\
David M. Mason\\
Department of Applied Economics and Statistics, University of Delaware \\
213 Townsend Hall, Newark, DE 19716, USA \\
e-mail: \texttt{davidm@udel.edu}}
\date{}
\maketitle

\begin{abstract} 
We define a time dependent empirical process based on $n$
independent fractional Brownian motions and describe strong approximations
to it by Gaussian processes. They lead to strong approximations and
functional laws of the iterated logarithm for the quantile or inverse of
this empirical process. They are obtained via time dependent Bahadur--Kiefer
representations.
\end{abstract}



\section{Introduction}

\label{s1} Swanson \cite{Swan07} using classical weak convergence theory
proved that an appropriately scaled median of $n$ independent Brownian
motions converges weakly to a mean zero Gaussian process. More recently
Kuelbs and Zinn \cite{KZ1}, \cite{KZ2} have obtained central limit theorems
for a time dependent quantile process based on $n$ independent copies of a
wide variety of random processes, which may be zero or perturbed to be not
zero with probability $1$ [w.p$.1$] at zero. These include certain
self-similar processes of which fractional Brownian motion is a special
case. Their approach is based on an extension of a result of Vervaat \cite%
{Vervaat} on the weak convergence of inverse processes in combination with
results from their deep study with Kurtz [Kurtz, Kuelbs and Zinn \cite{KKZ}]
of central limit theorems for time dependent empirical processes. \smallskip

\noindent We shall begin by defining a time dependent empirical process
based on $n$ independent fractional Brownian motions and describe a strong
approximations to it recently obtained by Kevei and Mason \cite{KM}. We
shall see that they lead to strong approximations and functional laws of the
iterated logarithm for the quantile or inverse of these empirical processes
and are obtained via time dependent Bahadur--Kiefer representations.

\subsection{Swanson (2007) result}

\noindent Our work is motivated by the following result of Swanson \cite%
{Swan07}.\smallskip

\noindent Let $\big\{ B_{j}^{(1/2)}\big\}_{j\geq 1}$ be a
sequence of i.i.d.~standard Brownian motions and for each $n\geq 1$ and $%
t\geq 0$ let $M_{n}\left( t\right) \mathbf{\ }$denote the median of$\
B_{1}^{\left( 1/2\right) }\left( t\right) ,\dots ,B_{n}^{\left( 1/2\right)
}\left( t\right) $. Swanson \cite{Swan07} using classical weak convergence
theory proved that $\sqrt{n}M_{n}\left( t\right) $ converges weakly to a
continuous centered Gaussian process $X$ on $\left[ 0,\infty \right) $ with
covariance function defined for $t_{1},t_{2}\in \left[ 0,\infty \right) $ by 
\begin{equation*}
E\left( X\left( t_{1}\right) X\left( t_{2}\right) \right) =\sqrt{t_{1}t_{2}}%
\sin ^{-1}\left( \frac{t_{1}\wedge t_{2}}{\sqrt{t_{1}t_{2}}}\right) .
\end{equation*}%
For a random particle motivation to look at such problems consult the
Introduction in \cite{Swan07}, where possible fractional Brownian motion
generalizations are hinted at.  \smallskip 

One of the aims of this paper is to place this result within the framework
of what has been long known about the usual empirical and quantile processes.

\subsection{Some classical quantile process lore}

To put our study into a broader context, we recall here some classical
quantile process lore. Let $X_{1},X_{2},\dots,$ be i.i.d.~$F$. For $\alpha
\in (0,1) $ define the inverse or quantile function $Q(\alpha) =\inf\left\{ 
x:F\left( x\right) \geq\alpha\right\} $ and the
empirical quantile function $Q_{n}\left( \alpha\right) =\inf\left\{
x:F_{n}(x) \geq\alpha\right\} $, where 
\begin{equation*}
F_{n} (x) =n^{-1}\sum_{j=1}^{n}1\left\{ X_{j}\leq x\right\}, \  
x\in\mathbb{R},
\end{equation*}
is the empirical distribution function based on $X_{1},\dots,X_{n}$.
\smallskip

We define the empirical process 
\begin{equation*}
v_{n}\left( x\right) :=\sqrt{n}\left\{ F_{n}\left( x\right) -F\left(
x\right) \right\} , \ x\in\mathbb{R},
\end{equation*}
and the quantile process%
\begin{equation*}
u_{n}\left( t\right) :=\sqrt{n}\left\{ Q_{n}\left( t\right) -Q\left(
t\right) \right\} ,\text{ }t\in\left( 0,1\right) .
\end{equation*}
For a real-valued function $\Upsilon$ defined on a set $S$ we
shall use the notation 
\begin{equation}
\left\Vert \Upsilon\right\Vert _{S}=\sup_{s\in S}\left\vert \Upsilon\left(
s\right) \right\vert . \label{sup}
\end{equation}
The empirical and quantile processes are closely connected to each other
through the following Bahadur--Kiefer representation:\smallskip

\noindent Let $X_{1},X_{2},\dots,$ be i.i.d.~$F$ on $\left[ 0,1\right] $,
where $F$ is twice differentiable on $\left( 0,1\right)$, $f(x) = F'(x)$, with 
\begin{equation*}
\inf_{x\in\left( 0,1\right) } f(x) >0 \text{ and}%
\sup_{x\in\left( 0,1\right) } \left\vert F^{\prime\prime}(x)\right\vert <\infty.
\end{equation*}
We have (Kiefer \cite{Kief}) the Bahadur--Kiefer representation 
\begin{equation}
\limsup_{n\rightarrow\infty}\frac{n^{1/4}\left\Vert v_{n}(Q)
+f(Q) u_{n}\right\Vert _{\left( 0,1\right) }}{\sqrt[4]{\log\log n%
}\sqrt{\log n}}=\frac{1}{\sqrt[4]{2}},\text{ a.s.}  \label{bk}
\end{equation}
The \textquotedblleft Bahadur\textquotedblright \ is in reference to the
original Bahadur \cite{B} paper, where a less precise version of (\ref{bk})
was first established. The function $f(Q)$ is called the
density quantile function. Deheuvels and Mason \cite{DeheuvelsMason}
developed a general approach to such theorems. For corresponding $L^p$ versions of such 
results we refer to Cs\"org\H{o} and Shi \cite{CsS}. 
\smallskip

\noindent Next using a strong approximation result of Koml\'{o}s, Major and
Tusn\'{a}dy \cite{KMT} one has on the same probability space an i.i.d.~%
$F$ sequence $X_{1},X_{2},\dots,$ and a sequence of i.i.d.~Brownian bridges $%
U_{1},U_{2},\dots,$ on $\left[ 0,1\right] $ such that 
\begin{equation}
\left\Vert v_{n}\left( Q\right) -\frac{\sum_{j=1}^{n}U_{j}}{\sqrt{n}}%
\right\Vert _{\left( 0,1\right) }=O\left( \frac{\left( \log n\right) ^{2}}{%
\sqrt{n}}\right) ,\text{ a.s.}  \label{kmt}
\end{equation}
Using (\ref{kmt}) it is easy see that under the conditions for which the
above the Bahadur--Kiefer representation (\ref{bk}) holds%
\begin{equation*}
\limsup_{n\rightarrow\infty}\frac{n^{1/4}\left\Vert \frac{\sum_{j=1}^{n}U_{j}%
}{\sqrt{n}}+f\left( Q\right) u_{n}\right\Vert _{\left( 0,1\right) }}{\sqrt[4]%
{\log\log n}\sqrt{\log n}}=\frac{1}{\sqrt[4]{2}},\text{ a.s.}  
\end{equation*}
Deheuvels \cite{Deheuvels} has shown that this rate of strong approximation
rate cannot be improved. \medskip

\noindent We shall develop analogues of these classical results for time
dependent empirical and quantile processes based on independent copies of
fractional Brownian motion. In particular, we shall extend the Swanson setup
to fractional Brownian motion, which will put his result in a broader
context.

\section{A time dependent empirical process}
\label{s2A}

In this section we recall some needed notation from \cite{KM}. Let 
$\left\{ B^{(H)}\right\} \cup \big\{ B_{j}^{(H)} \big\} _{j\geq1}$ be
a sequence of i.i.d. sample continuous fractional Brownian motions with
Hurst index $0<H<1$ defined on $[0,\infty)$. Note that $B^{(H)}$ is a
continuous mean zero Gaussian process on $[0,\infty)$ with covariance
function defined for any $s,t\in\lbrack0,\infty)$%
\begin{equation*}
E\left( B^{(H)}(s) B^{(H)}(t) \right) =\frac {1}{2}%
\left( \left\vert s\right\vert ^{2H}+\left\vert t\right\vert
^{2H}-\left\vert s-t\right\vert ^{2H}\right) .
\end{equation*}
By the L\'{e}vy modulus of continuity theorem for sample continuous
fractional Brownian motion $B^{(H)}$ with Hurst index $0<H<1$, (see
Corollary 1.1 of \cite{Wang}), we have for any $0<T<\infty$, w.p.~$1$, 
\begin{equation}
\sup_{0\leq s\leq t\leq T}\frac{\left\vert B^{(H)}\left( t\right)
-B^{(H)}\left( s\right) \right\vert }{f_{H}(t-s)}=:L<\infty,  \label{MC}
\end{equation}
where for $u\geq0$ 
\begin{equation}
f_{H}(u)=u^{H}\sqrt{1\vee\log u^{-1}}  \label{FH}
\end{equation}
and $a\vee b=\max\{a,b\}$. We shall take versions of $\left\{
B^{(H)}\right\} \cup \big\{ B_{j}^{(H)} \big\} _{j\geq1}$ such that (\ref%
{MC}) holds for all of their trajectories. \smallskip

For any $t\in\left[ 0,\infty\right) $ and $x\in\mathbb{R}$ let $F\left(
t,x\right) =P\left\{ B^{(H)}\left( t\right) \leq x\right\} .$ Note that 
\begin{equation}
F\left( t,x\right) =\Phi\left( x/t^{H}\right) ,  \label{Ft}
\end{equation}
where $\Phi\left( x\right) =P\left\{ Z\leq x\right\} ,$ with $Z$ being a
standard normal random variable. For any $n\geq1$ define the time dependent 
\textit{empirical distribution function} 
\begin{equation*}
F_{n}\left( t,x\right) =n^{-1}\sum_{j=1}^{n}1\left\{ B_{j}^{\left( H\right)
}\left( t\right) \leq x\right\} .
\end{equation*}

\noindent Applying Theorem 5 in \cite{KKZ} (also see their Remark 8) one can
show for any choice of $0<\gamma\leq1<T<\infty$ that the time dependent 
\textit{empirical process} indexed by $\left( t,x\right) \in\mathcal{T}%
\left( \gamma\right) $, 
\begin{equation*}
v_{n}\left( t,x\right) =\sqrt{n}\left\{ F_{n}\left( t,x\right) -F\left(
t,x\right) \right\} ,
\end{equation*}
where%
\begin{equation*}
\mathcal{T}\left( \gamma\right) :=\left[ \gamma,T\right] \times\mathbb{R},
\end{equation*}
converges weakly to a uniformly continuous centered Gaussian process $%
G\left( t,x\right) $ indexed by $\left( t,x\right) \in\mathcal{T}\left(
\gamma\right) $, whose trajectories are bounded, having covariance function%
\begin{equation}
\begin{split}
E\left( G\left( s,x\right) G\left( t,y\right) \right) 
=P\left\{ B^{\left( H\right) }\left( s\right) \leq x,B^{(H)}\left(
t\right) \leq y\right\} -P\left\{ B^{(H)}\left( s\right) \leq x\right\}
P\left\{ B^{(H)}\left( t\right) \leq y\right\} .
\end{split}
\label{EG}
\end{equation}

Here we restrict ourselves in stating this weak convergence result to
positive $\gamma$, since as pointed out in Section 8.1 of \cite{KKZ} the
empirical process $v_{n}(t,x) $ indexed by $\mathcal{T}\left(
0\right) :=\left[ 0,T\right] \times\mathbb{R}$ does not converge weakly to a
uniformly continuous centered Gaussian process indexed by $\left( t,x\right)
\in\mathcal{T}\left( 0\right) $, whose trajectories are bounded. In the
sequel, $G\left( t,x\right) $ denotes a centered Gaussian process on $%
\mathcal{T}\left( 0\right) $ with covariance (\ref{EG}) that is uniformly
continuous on $\mathcal{T}\left( \gamma\right) $ with bounded trajectories
for any $0<\gamma\leq1<T<\infty$.

We shall also be using the following empirical process indexed by function
notation. Let $X,X_{1},X_{2},\dots$, be i.i.d.~random variables from a
probability space $\left( \Omega,\mathcal{A},P\right) $ to a measurable
space $\left( S,\mathcal{S}\right) $. Consider an empirical process indexed
by a class $\mathcal{G}$ of bounded measurable real valued functions on $%
\left( S,\mathcal{S}\right) $ defined by 
\begin{equation*}
\alpha_{n}\left( \varphi\right) :=\sqrt{n}(P_{n}-P)\varphi=\frac{\sum
_{i=1}^{n}\varphi\left( X_{i}\right) -nE\varphi\left( X\right) }{\sqrt{n}}%
, \quad \varphi\in\mathcal{G},
\end{equation*}
where 
\begin{equation*}
P_{n}\left( \varphi\right) =n^{-1}\sum_{i=1}^{n}\varphi\left( X_{i}\right) 
\text{ and }P\left( \varphi\right) =E\varphi\left( X\right) \text{.}
\end{equation*}
Keeping this notation in mind, let $\mathcal{C}\left[ 0,T\right] $ be the
class of continuous functions $g$ on $\left[ 0,T\right] $ endowed with the
topology of uniform convergence. Define the
subclass of $\mathcal{C}\left[ 0,T\right] $ 
\begin{equation*}
\mathcal{C}_{\infty}:=\left\{ g:\ \sup\left\{ \frac{\left\vert g\left(
s\right) -g\left( t\right) \right\vert }{f_{H}(\left\vert s-t\right\vert )}%
,\ 0\leq s,t\leq T\right\} <\infty\right\} .  
\end{equation*}
Further, let $\mathcal{F}_{\left( \gamma,T\right) }$ be the class of
functions of $g\in$ $\mathcal{C}\left[ 0,T\right] \rightarrow\mathbb{R}$,
indexed by $\left( t,x\right) \in\mathcal{T}\left( \gamma\right) ,$ of the
form 
\begin{equation*}
h_{t,x}\left( g\right) =1\left\{ g\left( t\right) \leq x,g\in \mathcal{C}%
_{\infty}\right\} .
\end{equation*}
Here we permit $\gamma=0$. Since by (\ref{MC}) we can assume that each $%
B^{(H)},$ $B_{j}^{(H)}$, $j\geq1$, is in $\mathcal{C}_{\infty}$, we see that
for any $h_{t,x}\in\mathcal{F}_{\left( \gamma,T\right) }$, 
\begin{equation}
\alpha_{n}\left( h_{t,x}\right) =\frac{1}{\sqrt{n}}\sum_{i=1}^{n}\left(
1\left\{ B_{i}^{(H)}\left( t\right) \leq x\right\} -P\left\{ B^{(H)}\left(
t\right) \leq x\right\} \right) =v_{n}\left( t,x\right) .  \label{v}
\end{equation}
We shall be using the notation $\alpha_{n}\left( h_{t,x}\right) $ and $%
v_{n}\left( t,x\right) $ interchangeably. \smallskip

Let $\mathbb{G}_{\left( \gamma,T\right) }$ denote the mean zero Gaussian
process indexed by $\mathcal{F}_{\left( \gamma,T\right) }$, having
covariance function defined for $h_{s,x},h_{t,y}\in$ $\mathcal{F}_{\left(
\gamma,T\right) }$ 
\begin{equation*}
E\left( \mathbb{G}_{\left( \gamma,T\right) }\left( h_{s,x}\right) \mathbb{G}%
_{\left( \gamma,T\right) }\left( h_{t,y}\right) \right) =P\left\{ B^{\left(
H\right) }\left( s\right) \leq x,B^{(H)}\left( t\right) \leq y,B^{\left(
H\right) }\in\mathcal{C}_{\infty}\right\}
\end{equation*}%
\begin{equation*}
-P\left\{ B^{(H)}\left( s\right) \leq x,B^{(H)}\in\mathcal{C}_{\infty
}\right\} P\left\{ B^{(H)}\left( t\right) \leq y,B^{(H)}\in\mathcal{C}%
_{\infty}\right\} ,
\end{equation*}
which since $P\left\{ B^{(H)}\in\mathcal{C}_{\infty}\right\} =1$,%
\begin{equation*}
=E\left( G\left( s,x\right) G\left( t,y\right) \right) ,
\end{equation*}
i.e. $\mathbb{G}_{\left( \gamma,T\right) }\left( h_{t,x}\right) $ defines a
probabilistically equivalent version of the Gaussian process $G\left(
t,x\right) $ for $\left( t,x\right) \in\mathcal{T}\left( \gamma\right) $. We
shall say that a process $\widetilde{\mathcal{Y}}$ is a \textit{%
probabilistically equivalent version} of $\mathcal{Y}$ if $\widetilde{%
\mathcal{Y}}\overset{\mathrm{D}}{=}\mathcal{Y}$.

\subsection{The Kevei and Mason (2016) strong approximation results for $%
\protect\alpha_{n}$}

\label{main}For future reference we record here two strong approximations
for $\alpha_{n}$ that were recently established by Kevei and Mason \cite{KM}. In the 
results that follow 
\begin{equation}
\nu_{0}=2+\frac{2}{H}\quad\text{and }\ H_{0}=1+H.  \label{nuk}
\end{equation}

The main results in \cite{KM} are the following two strong approximation theorems.

\begin{theorem}(\cite{KM}) \label{th:1} 
For any $1 \geq  \gamma>0$, for all $1/\left( 2\tau
_{1}(0)\right) <\alpha<1/\tau_{1}(0)$ and $\xi>1$ there exist a $\rho\left(
\alpha,\xi\right) >0$, a sequence of i.i.d.~$B_{1}^{\left( H\right)
},B_{2}^{(H)},\ldots,$ and a sequence of independent copies $\mathbb{G}%
_{\left( \gamma,T\right) }^{\left( 1\right) },\mathbb{G}_{\left(
\gamma,T\right) }^{\left( 2\right) },\ldots$, of $\mathbb{G}_{\left(
\gamma,T\right) }$ sitting on the same probability space such that 
\begin{equation}
\max_{1\leq m\leq n}
\Big\Vert \sqrt{m}\alpha_m - \sum_{i=1}^{m} \mathbb{G}_{(\gamma,T)}^{(i)}
\Big\Vert_{\mathcal{F}_{(\gamma,T)}}
=O\left( n^{1/2-\tau\left( \alpha\right) }\left(\log n\right)^{\tau_{2}}\right),
\text{ a.s.,}  \label{tt1}
\end{equation}
where $\tau\left( \alpha\right) =\left( \alpha\tau_{1}(0)-1/2\right)
/(1+\alpha)>0$, $\tau_{1}(0)=1/\left( 2+5\nu_{0}\right) $, $\tau_{2}=(19H+25)/(24H+20)$ 
and $\nu_{0}$ is defined in (\ref{nuk}).
\end{theorem}

For any $\kappa>0$ let 
\begin{equation*}
\mathcal{G}\left( \kappa\right) =\left\{ t^{\kappa}h_{t,x}:\left( t,x\right)
\in\left[ 0,T\right] \times\mathbb{R}\right\} .  
\end{equation*}
For $g\in$ $\mathcal{G}\left( \kappa\right) $, with some abuse of notation,
we shall write%
\begin{equation*}
\mathbb{G}_{\left( 0,T\right) }\left( g\right) =t^{\kappa}\mathbb{G}_{\left(
0,T\right) }\left( h_{t,x}\right) .  
\end{equation*}
Also, in analogy with (\ref{sup}), in the following theorem, 
\begin{equation*}
\Big\Vert \sqrt{m}\alpha_{m} - \sum_{i=1}^{m}\mathbb{G}_{(0,T)}^{(i)}
\Big\Vert_{\mathcal{G}(\kappa)}:=
\sup \bigg\{ \Big\vert 
t^{\kappa}\alpha_{n}(h_{t,x}) - t^{\kappa} 
\sum_{i=1}^{m}\mathbb{G}_{(0,T)}^{(i)}\left( h_{t,x} \right) \Big\vert :
(t,x) \in [ 0,T] \times\mathbb{R}\bigg\} .
\end{equation*}

\begin{theorem}(\cite{KM}) \label{th:2} 
For any $\kappa >0$, for all $1/\left( 2\tau_{1}^{\prime
}\right) <\alpha <1/\tau _{1}^{\prime }$, and $\xi >1$ there exist a $\rho
^{\prime }\left( \alpha ,\xi \right) >0$, a sequence of i.i.d.~$%
B_{1}^{\left( H\right) },B_{2}^{(H)},\ldots ,$ and a sequence of independent
copies $\mathbb{G}_{\left( 0,T\right) }^{\left( 1\right) },\mathbb{G}%
_{\left( 0,T\right) }^{\left( 2\right) },$ $\dots ,$ of $\mathbb{G}_{\left(
0,T\right) }$ sitting on the same probability space such that 
\begin{equation}
\max_{1\leq m\leq n}
\Big\Vert \sqrt{m}\alpha_m - \sum_{i=1}^{m} \mathbb{G}_{(0,T)}^{(i)}
\Big\Vert_{\mathcal{G}(\kappa) }
= O \left( n^{1/2-\tau ^{\prime }\left( \alpha \right) }
\left(\log n\right)^{\tau_2}\right) \text{, a.s.,}  \label{t3}
\end{equation}%
where $\tau ^{\prime }\left( \alpha \right) =\left( \alpha \tau _{1}^{\prime
}-1/2\right) /(1+\alpha )>0$ and $\tau _{1}^{\prime }=\tau _{1}^{\prime
}(\kappa )=\kappa /(5H_{0}+\kappa (2+5\nu _{0})).$
\end{theorem}

Notice that (\ref{tt1}) and (\ref{t3}) trivially imply that for some $%
1/2>\xi>0$ 
\begin{equation*}
\max_{1\leq m\leq n} \Big\Vert \sqrt{m}\alpha_{m} - 
\sum_{i=1}^{m}\mathbb{G}_{(\gamma,T)}^{(i)}
\Big\Vert_{\mathcal{F}_{(\gamma,T)}}
=O\left( n^{-\xi}\right), \ \text{ a.s.,}
\end{equation*}
and 
\begin{equation*}
\max_{1\leq m\leq n}\Big\Vert \sqrt{m}\alpha_{m} - \sum_{i=1}^{m}
\mathbb{G}_{(0,T)}^{(i)} \Big\Vert_{\mathcal{G}(\kappa)}
=O\left( n^{-\xi}\right), \ \text{ a.s.}  
\end{equation*}

\subsection{Applications to LIL}

Kevei and Mason $\cite{KM}$ point out that the following compact law of the
iterated logarithm (LIL) for $\alpha_{n}$ follows from their Theorem \ref%
{th:1}, namely 
\begin{equation}
\left\{ \frac{\alpha_{n}\left( h_{t,x}\right) }{\sqrt{2\log\log n}}%
:h_{t,x}\in\mathcal{F}_{\left( \gamma,T\right) }\right\} =\left\{ \frac{%
v_{n}\left( t,x\right) }{\sqrt{2\log\log n}}:\left( t,x\right) \in\mathcal{T}%
\left( \gamma\right) \right\}  \label{LL2}
\end{equation}
is, w.p.~$1$, relatively compact in $\ell_{\infty}\left( \mathcal{F}_{\left(
\gamma,T\right) }\right) $ (the space of bounded functions $\Upsilon$ on $%
\mathcal{F}_{\left( \gamma,T\right) }$ equipped with supremum norm $%
\left\Vert \Upsilon\right\Vert _{\mathcal{F}_{\left( \gamma,T\right)
}}=\sup_{\varphi\in\mathcal{F}_{\left( \gamma,T\right) }}\left\vert
\Upsilon\left( \varphi\right) \right\vert $) and its limit set is the unit
ball of the reproducing kernel Hilbert space determined by the covariance
function $E (\mathbb{G}_{(\gamma,T)} (h_{s,x}) \mathbb{G}_{(\gamma,T)}(h_{t,y}) )$ 
$=E( G(s,x)  G(t,y))$. In particular we get that 
\begin{equation*}
\limsup_{n\rightarrow\infty}\frac{\left\Vert \alpha_{n}\right\Vert _{%
\mathcal{F}_{\left( \gamma,T\right) }}}{\sqrt{2\log\log n}}=\limsup
_{n\rightarrow\infty}\sup_{\left( t,x\right) \in\mathcal{T}\left(
\gamma\right) }\left\vert \frac{v_{n}\left( t,x\right) }{\sqrt{2\log\log n}}%
\right\vert =\sigma\left( \gamma,T\right) ,\text{ a.s.}  
\end{equation*}
where 
\begin{equation*}
\sigma^{2}\left( \gamma,T\right) =\sup\left\{ E\left( \mathbb{G}_{\left(
\gamma,T\right) }^{2}\left( h_{t,x}\right) \right) :h_{t,x}\in \mathcal{F}%
_{\left( \gamma,T\right) }\right\} =\frac{1}{4}.
\end{equation*}

Furthermore, they derive from their Theorem \ref{th:2} the following compact
LIL, for all $0<\kappa<\infty$, 
\begin{equation}
\left\{ \frac{t^{\kappa}\alpha_{n}\left( h_{t,x}\right) }{\sqrt{2\log\log n}}%
:h_{t,x}\in\mathcal{F}_{\left( 0,T\right) }\right\} =\left\{ \frac{%
t^{\kappa}v_{n}\left( t,x\right) }{\sqrt{2\log\log n}}:\left( t,x\right) \in%
\left[ 0,T\right] \times\mathbb{R}\right\}  \label{comL}
\end{equation}
is, w.p.~$1$, relatively compact in $\ell_{\infty}\left( \mathcal{G}\left(
\kappa\right) \right) $ and its limit set is the unit ball of the
reproducing kernel Hilbert space determined by the covariance function $%
E\left( s^{\kappa}t^{\kappa}\mathbb{G}_{\left( \gamma,T\right) }\left(
h_{s,x}\right) \mathbb{G}_{\left( \gamma,T\right) }\left( h_{t,y}\right)
\right) =E\left( s^{\kappa}t^{\kappa}G\left( s,x\right) G\left( t,y\right)
\right) .$ This implies that 
\begin{equation}
\limsup_{n\rightarrow\infty}\frac{\left\Vert \alpha_{n}\right\Vert _{%
\mathcal{G}\left( \kappa\right) }}{\sqrt{2\log\log n}}=\limsup
_{n\rightarrow\infty}\sup_{\left( t,x\right) \in\left[ 0,T\right] \times%
\mathbb{R}}\left\vert \frac{t^{\kappa}v_{n}\left( t,x\right) }{\sqrt{%
2\log\log n}}\right\vert =\sigma_{\kappa}\left( T\right) \text{, a.s.,}  \label{sig}
\end{equation}
where 
\begin{equation}
\sigma_{\kappa}^{2}\left( T\right) =\sup\left\{ E\left( \mathbb{G}_{\left(
0,T\right) }^{2}\left( t^{\kappa}h_{t,x}\right) \right) :t^{\kappa}h_{t,x}\in%
\mathcal{G}\left( \kappa\right) \right\} =\frac{T^{2\kappa}}{4}.
\label{sig2}
\end{equation}

\section{Bahadur--Kiefer representations and strong approximations for time
dependent quantile processes}
\label{s3}

\subsection{A time dependent quantile process}

For each $t\in\left( 0,\infty\right) $ and $\alpha\in\left( 0,1\right) $
define the time dependent \textit{inverse} or \textit{quantile function}%
\begin{equation*}
\tau_{\alpha}( t ) =\inf\left\{ x:F\left( t,x\right) \geq
\alpha\right\} ,
\end{equation*}
and the time dependent \textit{empirical inverse} or \textit{empirical
quantile function} 
\begin{equation} \label{eq:def-taun}
\tau_{\alpha}^{n} ( t ) =\inf\left\{ x:F_{n}\left( t,x\right)
\geq\alpha\right\},
\end{equation}
and the corresponding time dependent \textit{quantile process }%
\begin{equation*}
u_{n} ( t,\alpha ) :=\sqrt{n}\left( \tau_{\alpha}^{n}(t) -\tau_{\alpha}(t) \right) .
\end{equation*}
Notice that by (\ref{Ft}), for each fixed $t>0$, $F\left( t,x\right) $ has
density 
\begin{equation*}
f\left( t,x\right) =\frac{1}{t^{H}\sqrt{2\pi}}\exp\left( -\frac{x^{2}}{%
2t^{2H}}\right) \text{, }-\infty<x<\infty\text{.}
\end{equation*}
Further, for each $t\in\left( 0,\infty\right) $ and $\alpha\in\left(
0,1\right) $, $\tau_{\alpha}\left( t\right) $ is uniquely defined by 
\begin{equation} \label{eq:def-tau}
\tau_{\alpha}\left( t\right) =t^{H}z_{\alpha},\text{ where }P\left\{ Z\leq
z_{\alpha}\right\} =\alpha,
\end{equation}
which says that $f\left( t,\tau_{\alpha}\left( t\right) \right) =\frac {1}{%
t^{H}\sqrt{2\pi}}\exp\left( -\frac{z_{\alpha}^{2}}{2}\right) .$

\subsection{Our results for time dependent quantile processes}

We shall prove the following uniform time dependent Bahadur--Kiefer
representations for the quantile process $u_{n}\left( t,\alpha\right) $. We
shall see that one easily infers from them LIL and strong approximations for
such processes.

Introduce the condition on a sequence of constants $0<\gamma_{n}\leq1$ 
\begin{equation}
\infty>-\frac{\log\gamma_{n}}{\log n}\rightarrow\eta,\ \text{ as }%
n\rightarrow\infty.  \label{eta}
\end{equation}

\begin{theorem}
\label{th:3} Whenever $0<\gamma=\gamma_{n}\leq1$ satisfies (\ref{eta}) for
some $0 \leq \eta < 1/(2 H)$, then for any $0<\rho<1/2$ and $T > 1$
\begin{equation} \label{bah}
\begin{split}
& \sup_{\left( t,\alpha\right) \in\left[ \gamma_{n},T\right] \times\left[
\rho,1-\rho\right] }\left\vert v_n (t,\tau_\alpha (t) ) +f ( t,\tau_{\alpha}(t) ) 
u_{n}(t,\alpha) \right\vert  \\
& =O\left( n^{-1/4}\gamma_{n}^{-H/2}\left( \log\log n\right)
^{1/4}\left( \log n\right) ^{1/2}\right) ,\text{ a.s.}  
\end{split}
\end{equation}
\end{theorem}

\begin{remark}
It is noteworthy here to point out that when $\gamma_{n}=\gamma$ is
constant, the rate in (\ref{bah}) corresponds to the known exact rate 
in (\ref{bk}) in the classic uniform Bahadur--Kiefer representation of sample quantiles. 
Refer to Deheuvels and Mason \cite{DeheuvelsMason} for more results in this direction.
\end{remark}

\begin{remark}
Note that in (\ref{bah}) smaller $\delta$ implies better rate, so it is enough to 
prove the statement for $0 < \delta$ small enough. Furthermore, if $\gamma_n \to 0$ 
and $\delta$ large it can happen that the rate in (\ref{bah}) tends to infinity. Since 
(\ref{eta}) holds with $\eta< 1/(2H)$, $\gamma_n$ cannot tend to zero too fast. In the 
borderline case (which is not allowed in the statement of the theorem) $\gamma_n = 
n^{-1/(2H)}$, the rate would go to infinity.
\end{remark}

\begin{remark}
Let $\ell_{\infty}\left( \left[ \gamma,T\right] \times\left[ \rho ,1-\rho%
\right] \right) \mathcal{\ }$denote the class of bounded functions on $\left[
\gamma,T\right] \times\left[ \rho,1-\rho\right] $. Notice when $%
0<\gamma\leq1 $ is fixed, we immediately get from (\ref{LL2}) and (\ref{bah}) that 
\begin{equation*}
\left\{ \frac{f( t,\tau_{\alpha}(t) )\, u_{n}(t,\alpha) }{\sqrt{2\log\log n}}:\left( 
t,\alpha\right) \in\left[
\gamma,T\right] \times\left[ \rho,1-\rho\right] \right\}
\end{equation*}
is, w.p.~$1$, relatively compact in $\ell_{\infty}\left( \left[ \gamma,T%
\right] \times\left[ \rho,1-\rho\right] \right) $ and its limit set is the
unit ball of the reproducing kernel Hilbert space determined by the
covariance function defined for $(t_{1},\alpha_{1})$, $(t_{2},\alpha_{2})$
$\in [\gamma,T] \times [\rho ,1-\rho]$ by 
\begin{align*}
K\left( \left( t_{1},\alpha_{1}\right) ,\left( t_{2},\alpha_{2}\right)
\right) & =E\left( G\left( t_{1},\tau_{\alpha_{1}}\left( t_{1}\right)
\right) G\left( t_{2},\tau_{\alpha_{2}}\left( t_{2}\right) \right) \right) \\
& =P\left\{ B^{\left( H\right) }\left( t_{1}\right) \leq
t_{1}^{H}z_{\alpha_{1}},B^{\left( H\right) }\left( t_{2}\right) \leq
t_{2}^{H}z_{\alpha_{2}}\right\} -\alpha_{1}\alpha_{2}.
\end{align*}
Also we get when $0<\gamma\leq1$ is fixed the following strong
approximation, namely on the probability space of Theorem \ref{th:1},%
\begin{equation*}
\sup_{\left( t,\alpha\right) \in\left[ \gamma,T\right] \times\left[
\rho,1-\rho\right] }\left\vert \sqrt{n} \, f( t,\tau_{\alpha}(t))\, u_{n}(t,\alpha) 
+\sum_{i=1}^{n}G_{i}(t,\tau_{\alpha}(t) ) \right\vert =O\left( 
n^{1/2-\tau
\left( \alpha\right) }\left( \log n\right) ^{\tau_{2}}\right) \text{, a.s.,}  
\end{equation*}
where $G_{i}(t,\tau_{\alpha}(t)) = 
\mathbb{G}_{(\gamma,T)}^{(i)} (h_{t,\tau_{\alpha}(t)})$. This follows from Theorems 
\ref{th:1} and \ref{th:3} by noting $\tau(\alpha)<1/4$.
\end{remark}

\begin{corollary} \label{cor1} 
For any $0<\rho<1/2$, $T>1$ and $\delta > 0$ we have
\begin{equation}
\sup_{\left( t,\alpha\right) \in\left[ 0,T\right] \times\left[ \rho,1-\rho%
\right] }\left\vert t^{H}v_{n}\left( t,\tau_{\alpha}\left( t\right) \right) +%
\frac{\exp\left( -\frac{z_{\alpha}^{2}}{2}\right) }{\sqrt{2\pi}}u_{n}\left(
t,\alpha\right) \right\vert =O\left( n^{-1/6 +\delta} \right) ,\text{ a.s.}
\label{c1}
\end{equation}
\end{corollary}

\begin{remark}
Let $\ell_{\infty}\left( \left[ 0,T\right] \times\left[ \rho ,1-\rho\right]
\right) \mathcal{\ }$denote the class of bounded functions on $\left[ 0,T%
\right] \times\left[ \rho,1-\rho\right] $. Observe that (\ref{c1}) combined
with the compact LIL pointed out in (\ref{comL}), immediately imply that 
\begin{equation*}
\left\{ \frac{\exp\left( -\frac{z_{\alpha}^{2}}{2}\right) u_{n}\left(
t,\alpha\right) }{\sqrt{2\pi}\sqrt{2\log\log n}}:\left( t,\alpha\right) \in%
\left[ 0,T\right] \times\left[ \rho,1-\rho\right] \right\}
\end{equation*}
is, w.p.~$1$, relatively compact in $\ell_{\infty}\left( \left[ 0,T\right]
\times\left[ \rho,1-\rho\right] \right) $ and its limit set is the unit ball
the reproducing kernel Hilbert space determined by the covariance function
defined for $(t_{1},\alpha_{1}),(t_{2},\alpha_{2})\in\lbrack
0,T]\times\lbrack\rho,1-\rho]$ by 
\begin{align*}
K\left( \left( t_{1},\alpha_{1}\right) ,\left( t_{2},\alpha_{2}\right)
\right) & =t_{1}^{H}t_{2}^{H}E\left( G( t_{1},\tau_{\alpha_{1}}(t_{1})) 
G(t_{2},\tau_{\alpha_{2}}(t_{2})) \right) \\
& =t_{1}^{H}t_{2}^{H}\left( P\left\{ B^{\left( H\right) }( t_{1})
\leq t_{1}^{H}z_{\alpha_{1}},B^{\left( H\right) }(t_{2}) \leq
t_{2}^{H}z_{\alpha_{2}}\right\} -\alpha_{1}\alpha _{2}\right) .
\end{align*}
We also get the following strong approximation, namely on the probability
space of Theorem \ref{th:2} with $\kappa=H$, for some $1/2>\xi>0$ 
\begin{equation}
\sup_{\left( t,\alpha\right) \in\left[ 0,T\right] \times\left[ \rho,1-\rho%
\right] }\left\vert \frac{\exp\left( -\frac{z_{\alpha}^{2}}{2}\right)
u_{n}( t,\alpha) }{\sqrt{2\pi}}+\frac{1}{\sqrt{n}}%
\sum_{i=1}^{n}t^{H}G_{i} ( t,\tau_{\alpha}(t))
\right\vert =O\left( n^{-\xi}\right), \text{ a.s.,}  \label{sau}
\end{equation}
where $G_{i}(t,\tau_{\alpha}(t))=\mathbb{G}_{(0,T)}^{(i)}\left(
h_{t,\tau_{\alpha}(t)}\right) $. This follows from (\ref{t3}), noting that $%
\tau^{\prime}(\alpha)>0$, combined with (\ref{c1}).
\end{remark}

\begin{remark}
Let $\ell_{\infty}\left( \left[ 0,T\right] \right) \mathcal{\ }$denote the
class of bounded functions on $\left[ 0,T\right] $. Applying the compact LIL
pointed out in the previous remark with $H=1/2,$ to the median process
considered by Swanson \cite{Swan07}, i.e. 
\begin{equation*}
\sqrt{n}M_{n}\left( t\right) =u_{n}\left( t,1/2\right) =\sqrt{n}\tau
_{1/2}^{n}\left( t\right) ,\text{ }t\geq0\text{,}
\end{equation*}
we get for any $T>0$, that 
\begin{equation*}
\left\{ \frac{\sqrt{n}M_{n}\left( t\right) }{\sqrt{2\log\log n}}:t\in\left[
0,T\right] \right\}
\end{equation*}
is, w.p.~$1$, relatively compact in $\ell_{\infty}\left( \left[ 0,T\right]
\right) ,$ and its limit set is the unit ball of the reproducing kernel
Hilbert space determined by the covariance function defined for $%
t_{1},t_{2}\in\left[ 0,T\right] $ 
\begin{align*}
2\pi K\left( t_{1},t_{2}\right) & =2\pi\sqrt{t_{1}t_{2}}E\left( G\left(
t_{1},0\right) G\left( t_{2},0\right) \right) \\
& =2\pi\sqrt{t_{1}t_{2}}\left( P\left\{ B^{\left( 1/2\right) }\left(
t_{1}\right) \leq0,B^{\left( 1/2\right) }\left( t_{2}\right) \leq0\right\}
-1/4\right) ,
\end{align*}
which equals 
\begin{equation}
\sqrt{t_{1}t_{2}}\sin^{-1}\left( \frac{t_{1}\wedge t_{2}}{\sqrt{t_{1}t_{2}}}%
\right) .  \label{xco}
\end{equation}
In particular we get%
\begin{equation*}
\limsup_{n\rightarrow\infty}\frac{\left\Vert \sqrt{n}M_{n}\right\Vert _{%
\left[ 0,T\right] }}{\sqrt{2\log\log n}}=\sqrt{T\sin^{-1}\left( 1\right) }=%
\sqrt{T\pi/2},\text{ a.s.}
\end{equation*}
Moreover, since a mean zero Gaussian process $X(t)$, $t\geq 0$,
with covariance function (\ref{xco}) is equal in distribution to $-\sqrt{2\pi 
t}G(t,0)$, $t\geq0$, we see from (\ref{sau}) that there
exist a sequence $B_{1}^{\left( 1/2\right) },B_{2}^{\left( 1/2\right)
},\ldots,$ i.i.d.~$B^{\left( 1/2\right) }$ and a sequence of processes $%
X^{\left( 1\right) },X^{\left( 2\right) },\ldots$, i.i.d.~$X$ sitting on the
same probability space such that, a.s. 
\begin{equation*}
\left\Vert \sqrt{n}M_{n}-\frac{1}{\sqrt{n}}\sum_{i=1}^{n}X^{\left( i\right)
}\right\Vert _{\left[ 0,T\right] }=o\left( 1\right) .
\end{equation*}
Of course, this implies the Swanson result that $\sqrt{n}M_{n}$ converges
weakly on $[ 0,T] $ to the process $X$.
\end{remark}

\section{Proofs of Theorem \protect\ref{th:3} and Corollary \protect\ref{cor1}}
\label{s5}

To ease the notation we suppress the upper index from the fractional Brownian motions, 
that is, in the following $B, B_1, B_2, \ldots$ are i.i.d.~fractional Brownian motions 
with Hurst index $H$.

\subsection{Proof of Theorem \protect\ref{th:3}}

Before we can prove Theorem \ref{th:3} we must first gather together some
facts about $\tau_{\alpha}^{n} (t)$, defined in (\ref{eq:def-taun}).

\begin{proposition} \label{prop:3} With probability $1$ for any choice of $0<\rho<1/2$ 
uniformly in $t>0$, $n\geq1$ and $0<\rho\leq\alpha\leq1-\rho$ 
\begin{equation*}
0\leq F_{n}\left( t,\tau_{\alpha}^{n}\left( t\right) \right) -\alpha\leq m/n,
\end{equation*}
where $m=2(\left\lceil 2/H\right\rceil +1)$.
\end{proposition}

\noindent\textit{Proof} We require a lemma.

\begin{lemma} \label{lemma:1} Let $B_{j}$, $j=1,\dots,n,$ be i.i.d.~fractional Brownian 
motions on $\left[ 0,\infty\right) $ with Hurst index $0<H<1$, where $n\geq2\left\lceil 
2/H\right\rceil +2$, then w.p.~zero does there exist a subset $\left\{ 
i_{1},\dots,i_{m}\right\} \subset\left\{
1,\dots,n\right\} $, where $m=2\left\lceil 2/H\right\rceil +2$, such that for
some $t>0$%
\begin{equation*}
B_{i_{1}}( t) =\dots=B_{i_{m}}(t).
\end{equation*}
\end{lemma}

\noindent\textit{Proof} If such a subset exists then the paths of the
independent fractional Brownian motions in $\mathbb{R}^{k}$ with $2k=m$, 
\begin{equation}
X^{1}=\left( B_{i_{1}},\ldots, B_{i_{k}}\right) \text{ and }
X^{2}=\left( B_{i_{k+1}},\ldots, B_{i_{2k}} \right)  \label{x1x2}
\end{equation}
would have non-empty intersection except at $0$, which, since $k>2/H$,
contradicts the following special case of Theorem 3.2 in Xiao \cite{Xiao}:

\begin{theorem*}\textbf{(Xiao)}
Let $X^{1}(t)$, $t\geq0$, and 
$X^{2}(t)$, $t\geq0$, be two independent fractional Brownian motions in $\mathbb{R}^{d}$ 
with index $0<H<1$. If $2/H\leq d$, then w.p.~$1$, 
\begin{equation*}
X^{1}\left( \left[ 0,\infty\right) \right) \cap X^{2}\left( \left(
0,\infty\right) \right) =\emptyset.
\end{equation*}
\end{theorem*}

We apply this result with $X^{1}$ and $X^{2}$ as in (\ref{x1x2}). \hfill$\square$ 
\smallskip

Returning to the proof of Proposition \ref{prop:3}, choose $%
n\geq2\left\lceil 2/H\right\rceil +2$ and for any choice of $t>0$ let $%
B_{\left( 1\right) }(t) \leq\dots\leq
B_{\left( n\right) }(t) $ denote the order statistics of $B_{1}( t),\ldots, B_{n}(t)$. We 
see that for any $\alpha\in\left( 0,1\right)$, 
\begin{equation*}
F_{n}\left( t,B_{^{\left( \left\lceil \alpha n\right\rceil \right)
}} ( t ) \right) \geq\left\lceil \alpha
n\right\rceil /n\geq\alpha
\end{equation*}
and%
\begin{equation*}
F_{n}\left( t,B_{^{\left( \left\lceil \alpha n\right\rceil \right)
}} (t) -\right) \leq\left( \left\lceil \alpha
n\right\rceil -1\right) /n<\alpha.
\end{equation*}
Thus 
\begin{equation*}
\tau_{\alpha}^{n}\left( t\right) =\inf\left\{ x:F_{n}\left( t,x\right)
\geq\alpha\right\} =B_{^{\left( \left\lceil \alpha n\right\rceil \right)}} (t) .
\end{equation*}
Since by the above lemma, w.p.~$1$, for all $t>0$ 
\begin{equation*}
\sum_{j=1}^{n}1\left\{ B_{j}\left( t\right) =B_{^{\left(
\left\lceil \alpha n\right\rceil \right) }}\left(
t\right) \right\} <m=2\left\lceil 2/H\right\rceil +2,
\end{equation*}
we see that 
\begin{equation*}
\alpha\leq\left\lceil \alpha n\right\rceil /n\leq F_{n}\left( t,\tau_{\alpha
}^{n}\left( t\right) \right) \leq\left( \left\lceil \alpha n\right\rceil
+m-1\right) /n\leq\alpha+m/n.
\end{equation*}
Thus w.p.~$1$ for any choice of $0<\rho<1/2$ uniformly in $t>0$, $%
n\geq2\left\lceil 2/H\right\rceil +2$ and $0<\rho\leq\alpha\leq1-\rho$ 
\begin{equation*}
0\leq F_{n}\left( t,\tau_{\alpha}^{n}\left( t\right) \right) -\alpha\leq m/n%
\text{.}
\end{equation*}
Note that this bound is trivially true for $1\leq n<2\left\lceil
2/H\right\rceil +2$. \hfill$\square$

\begin{proposition} \label{prop:4} For any $H\geq \delta>0$ and $\rho\in(0,1/2)$ there is 
a $D_{0}= D_0(\rho, T)>0$ (depending only on $\rho$ and $T$) such that,
w.p.~$1$ there is an $n_{0}=n_{0}(\delta)$, such that for all $n>n_{0}$, uniformly in 
$(\alpha,t)\in\lbrack\rho ,1-\rho]\times(a_{n}(\delta),T],$ 
\begin{equation*}
\left\vert \tau_{\alpha}\left( t\right) -\tau_{\alpha}^{n}\left( t\right) \right\vert 
\leq \frac{t^{H-\delta}D_{0}\sqrt{\log\log n}}{\sqrt{n}},
\end{equation*}
with 
\begin{equation}
a_{n}=a_{n}(\delta)=C\left( \frac{\log\log n}{n}\right) ^{1/(2\delta)},
\label{eq:a_n}
\end{equation}
where $C=C(\delta, \rho, T)$ depends only on  $\delta, \rho$ and $T$.
\end{proposition}

\noindent\textit{Proof} By Proposition \ref{prop:3}, w.p.~1, 
\begin{equation}
\sup_{\left( \alpha,t\right) \in\left[ \rho,1-\rho\right] \times\left( 0,T%
\right] }\left\vert F_{n}\left( t,\tau_{\alpha}^{n}\left( t\right) \right)
-\alpha\right\vert \leq m/n.  \label{Cn}
\end{equation}
We see by (\ref{sig}) that for any $H \geq \delta>0$ w.p.~1 there is an $n_{0}$, such
that for all $n>n_{0}$ 
\begin{equation*}
\sup_{\left( \alpha,t\right) \in\left[ \rho,1-\rho\right] \times\left( 0,T%
\right] }t^{\delta}\left\vert F_{n}\left( t,\tau_{\alpha}^{n}\left( t\right)
\right) -F\left( t,\tau_{\alpha}^{n}\left( t\right) \right) \right\vert \leq%
\frac{2\sigma_{\delta}(T)\sqrt{\log\log n}}{\sqrt{n}},
\end{equation*}
where, as in (\ref{sig2}), $\sigma_{\delta}^{2}\left( T\right) =\frac{%
T^{2\delta}}{4}\leq\frac{T^{2}}{4}.$ Thus by (\ref{Cn}) and noting that $%
F\left( t,\tau_{\alpha}\left( t\right) \right) =\alpha$ we have w.p.~1 for all
large enough $n$ 
\begin{equation}
\sup_{(\alpha,t)\in\lbrack\rho,1-\rho]\times(0,T]}t^{\delta}\left\vert
F\left( t,\tau_{\alpha}\left( t\right) \right) -F\left( t,\tau_{\alpha
}^{n}\left( t\right) \right) \right\vert \leq\frac{2T\sqrt{\log\log n}}{%
\sqrt{n}}.  \label{C3}
\end{equation}
Recall the notation in (\ref{eq:def-tau}).
Notice that whenever $t^{H}x-\tau_{\alpha}\left( t\right) >t^{H}/8$, for
some $t>0$ and $\alpha\in\left[ \rho,1-\rho\right] ,$ 
\begin{equation*}
\begin{split}
\left\vert F\left( t,\tau_{\alpha}\left( t\right) \right) -F\left(
t,t^{H}x\right) \right\vert 
& =\int_{\tau_{\alpha}\left( t\right) }^{t^{H}x}%
\frac{1}{t^{H}\sqrt{2\pi}}\exp\left( -\frac{y^{2}}{2t^{2H}}\right) dy \\
& =\int_{z_{\alpha}}^{x}\frac{1}{\sqrt{2\pi}}\exp\left( -\frac{u^{2}}{2}%
\right) du \\
& >\int_{z_{\alpha}}^{z_{\alpha}+1/8}\frac{1}{\sqrt{2\pi}}\exp\left(
-\frac{u^{2}}{2}\right) du\geq d_{1}>0,
\end{split}
\end{equation*}
where 
\begin{equation*}
d_{1}=\inf\left\{ \int_{z_{\alpha}}^{z_{\alpha}+1/8}\frac{1}{\sqrt{2\pi}}%
\exp\left( -\frac{u^{2}}{2}\right) du:\alpha\in\left[ \rho,1-\rho\right]
\right\} .
\end{equation*}
Similarly, whenever $\tau_{\alpha}\left( t\right) -t^{H}x>t^{H}/8$ for some $%
t>0$ and $\alpha\in\left[ \rho,1-\rho\right] ,$ 
\begin{equation*}
\left\vert F\left( t,\tau_{\alpha}\left( t\right) \right) -F\left(
t,t^{H}x\right) \right\vert \geq d_{1}.
\end{equation*}

We have shown that whenever $|t^{H}x-\tau_{\alpha}\left( t\right) |>t^{H}/8$, for some 
$t>0$, and $\alpha\in\left[ \rho,1-\rho\right] $, then 
\begin{equation*}
|F(t,\tau_{\alpha}(t))-F(t,t^{H}x)|>d_{1}>0.
\end{equation*}
Choose $C(\delta, \rho, T) = \left( 2T/d_{1}\right)^{1/\delta}$ in (\ref{eq:a_n}).
Then
\begin{equation}
\frac{2T\sqrt{\log\log n}}{\sqrt{n}}a_{n}^{-\delta}=\frac{2T}{C^{\delta}}%
= d_{1}.  \label{d1}
\end{equation}
Now, (\ref{C3}) implies that w.p.~$1$ for all large $n$
we have $|\tau_{\alpha}(t)-\tau_{\alpha}^{n}(t)|\leq t^{H}/8$, 
whenever $t>a_{n}$, which together
with $\alpha \in [ \rho,1-\rho]$ implies that 
\begin{equation}
\tau_{\alpha}(t),\tau_{\alpha}^{n}(t)\in t^{H}[z_{\rho}-1/8,z_{1-\rho
}+1/8]=:t^{H}\left[ a,b\right] .  \label{ab}
\end{equation}
We get for $t>a_{n}$ 
\begin{equation*}
|F(t,\tau_{\alpha}(t))-F(t,\tau_{\alpha}^{n}(t))|=\left\vert \Phi(\tau
_{\alpha}(t)t^{-H})-\Phi(\tau_{\alpha}^{n}(t)t^{-H})\right\vert
=t^{-H}|\tau_{\alpha}(t)-\tau_{\alpha}^{n}(t)|\varphi(\xi),
\end{equation*}
where $\xi\in\lbrack z_{\rho}-1/8,z_{1-\rho}+1/8]$, $\varphi$ is the
standard normal density and 
\begin{equation*}
\varphi(\xi)\geq\min_{y\in\lbrack a,b]}\varphi(y)=:d_{2}>0.
\end{equation*}
Therefore by (\ref{C3}), w.p.~$1$, for all large $n$, for $t>a_{n}$ and $%
\alpha\in\lbrack\rho,1-\rho]$ 
\begin{equation*}
|\tau_{\alpha}(t)-\tau_{\alpha}^{n}(t)|\leq\frac{2T}{d_{2}}\frac{t^{H-\delta
}\sqrt{\log\log n}}{\sqrt{n}},
\end{equation*}
so the statement is proved, with $D_{0}=2T/d_{2}$. \hfill$\square$
\medskip

For future reference we point out here that for any $a_{n}\left(
\delta\right)  $ as in (\ref{eq:a_n}) and $1\geq$ $\gamma_{n}>0$ satisfying
(\ref{eta}) for some $\eta<\frac{1}{2H}$
\[
\lim_{n\rightarrow\infty}\frac{-\log a_{n}\left(  \delta\right)  }{\log
n}=\frac{1}{2\delta} \geq \frac{1}{2H}>\lim_{n\rightarrow\infty}\frac{-\log
\gamma_{n}}{\log n}=\eta.
\]
Thus for all $n$ sufficiently large
\begin{equation}
a_{n}\left(  \delta\right)  <\gamma_{n}.\label{less}%
\end{equation}

Note that 
\begin{equation}
v_{n}\left( t,\tau_{\alpha}^{n}\left( t\right) \right) -\sqrt{n}\left\{
\alpha-F\left( t,\tau_{\alpha}^{n}\left( t\right) \right) \right\} =\sqrt{n}%
\left( F_{n}\left( t,\tau_{\alpha}^{n}\left( t\right) \right) -\alpha\right)
=:\Delta_{n}\left( t,\alpha\right) ,  \label{va}
\end{equation}
for which by Proposition \ref{prop:3} we have 
\begin{equation} \label{eq:Deltain1}
\left\vert \Delta_{n}\left(
t,\alpha\right) \right\vert \leq \frac{m}{\sqrt{n}},
\text{ uniformly in } t>0, \ 0<\rho \leq\alpha\leq1-\rho
\text{ and  } n\geq1.
\end{equation}
Rewriting (\ref{va}) as 
\begin{equation*}
v_{n}(t,\tau_{\alpha}^{n}(t))=-\sqrt{n}\left\{
F(t,\tau_{\alpha}^{n}(t))-\alpha\right\} +\Delta_{n}(t,\alpha),
\end{equation*}
we get using a Taylor expansion applied to $F(t,\tau_{\alpha}^{n}(t))-\alpha 
$, 
\begin{equation} \label{eq:v1}
v_{n}(t,\tau_{\alpha}^{n}(t))
=-\sqrt{n}f(t,\tau_{\alpha}(t))\left( \tau_{\alpha}^{n}(t)-\tau_{\alpha
}(t)\right) -\frac{1}{2}\sqrt{n}f^{\prime}\left( t,\theta_{\alpha}^{n}\left(
t\right) \right) \left( \tau_{\alpha}^{n}\left( t\right)
-\tau_{\alpha}\left( t\right) \right) ^{2}+\Delta_{n}(t,\alpha),
\end{equation}
where $\theta_{\alpha}^{n}\left( t\right) $ is between $\tau_{\alpha}(t)$ and 
$\tau_{\alpha}^{n}\left( t\right) $ and $f^{\prime}\left(
t,x\right) =\partial f(t,x) /\partial x$. \smallskip Write%
\begin{equation*}
\sqrt{n}f^{\prime}\left( t,\theta_{\alpha}^{n}\left( t\right) \right) \left(
\tau_{\alpha}^{n}\left( t\right) -\tau_{\alpha}\left( t\right) \right) ^{2}=%
\sqrt{n}t^{2H}f^{\prime}\left( t,\theta_{\alpha}^{n}\left( t\right) \right)
t^{-2H}\left( \tau_{\alpha}^{n}\left( t\right) -\tau_{\alpha}\left( t\right)
\right) ^{2}.
\end{equation*}
Observe that by (\ref{less}) with $\left[ a,b\right] $ as given in (\ref{ab}), w.p.~1, for 
all large $n$ 
\begin{equation*} 
\begin{split}
 \sup\left\{ \left\vert t^{2H}f^{\prime}\left( t,\theta_{\alpha}^{n}\left(
t\right) \right) \right\vert :\left( \alpha,t\right) \in\left[ \rho,1-\rho%
\right] \times\left[ \gamma_{n},T\right] \right\} 
& \leq\sup_{t\in\left( 0,T\right] }\sup\left\{ t^{2H}\left\vert f^{\prime
}\left( t,x\right) \right\vert :x\in t^{H}\left[ a,b\right] \right\} \\
& = \sup\left\{ \left\vert f^{\prime}\left( 1,x\right) \right\vert :x\in 
\left[ a,b\right] \right\} <\infty.
\end{split}
\end{equation*}
Further by (\ref{less}), we can apply Proposition \ref{prop:4} with $\delta = H/4$ to get,
w.p.~1, 
\begin{equation*} 
\sup_{\left( \alpha,t\right) \in\left[ \rho,1-\rho\right] \times\left[
\gamma_{n},T\right] }t^{-2H}\left( \tau_{\alpha}^{n}\left(
t\right) -\tau_{\alpha}\left( t\right) \right) ^{2}=O\left( \frac {%
\gamma_{n}^{-H/2}\log\log n}{n}\right) .
\end{equation*}
Therefore, substituting back into (\ref{eq:v1}) from the definition of $u_n$ and from
(\ref{eq:Deltain1}) we see that w.p.~1,  
\begin{equation}
\sup_{\left( \alpha,t\right) \in\left[ \rho,1-\rho\right] \times\left[
\gamma_{n},T\right] }\left\vert v_{n}\left( t,\tau_{\alpha}^{n}\left(
t\right) \right) +f\left( t,\tau_{\alpha}\left( t\right) \right) u_{n}\left(
t,\alpha\right) \right\vert =O\left( \frac{\gamma_{n}^{-H/2}\log\log n}{%
\sqrt{n}}\right) .  \label{BHD}
\end{equation}
\smallskip 

Our next goal is to control the size of $v_{n}\left( t,\tau_{\alpha}\left(
t\right) \right) -v_{n}\left( t,\tau_{\alpha}^{n}\left( t\right) \right) $
uniformly in $\left( \alpha,t\right) \in\left[ \rho,1-\rho\right] \times%
\left[ \gamma_{n},T\right] $ for appropriate $0<\gamma_{n}\leq1$. For this
purpose we need to introduce some more notation.\smallskip

Recall notation (\ref{FH}). For any $K\geq1$ denote the class of real-valued
functions on $\left[ 0,T\right] $, 
\begin{equation*}
\mathcal{C}\left( K\right) =\left\{ g:\ \left\vert g(s)-g(t)\right\vert \leq
Kf_{H}(|s-t|),0\leq s,t\leq T\right\} .
\end{equation*}
One readily checks that $\mathcal{C}\left( K\right) $ is closed in $\mathcal{%
C}\left[ 0,T\right] $. The following class of functions $\mathcal{C}\left[
0,T\right] \rightarrow\mathbb{R}$ will play an essential role in our proof: 
\begin{equation*}
\mathcal{F}\left( K,\gamma\right) :=\left\{ h_{t,x}^{\left( K\right) }\left(
g\right) =1\left\{ g\left( t\right) \leq x,g\in\mathcal{C}\left( K\right)
\right\} :\left( t,x\right) \in\mathcal{T}\left( \gamma\right) \right\} .
\end{equation*}
For any $c>0$, $n>e$ and $1<T$ denote the class of real-valued
functions on $\left[ 0,T\right] $, 
\begin{equation}
\mathcal{C}_{n}:=\mathcal{C}(\sqrt{c\log n})=\left\{ g:\,\left\vert g\left(
s\right) -g\left( t\right) \right\vert \leq\sqrt{c\log n}f_{H}(\left\vert
s-t\right\vert ),\,0\leq s,t\leq T\right\} .  \label{cn}
\end{equation}
Define the class of functions $\mathcal{C}\left[ 0,T\right] \rightarrow 
\mathbb{R}$ indexed by $\left[ \gamma_{n},T\right] \times\mathbb{R=}\mathcal{%
T}\left( \gamma_{n}\right) $ 
\begin{equation*}
\mathcal{F}_n=\left\{ h_{t,x}^{\left( 
\sqrt{c\log n}\right) }\left( g\right) =1\left\{ g\left( t\right) \leq x,g\in%
\mathcal{C}_{n}\right\} :\left( t,x\right) \in\mathcal{T}\left(
\gamma_{n}\right) \right\} . 
\end{equation*}
To simplify our previous notation we shall write here 
\[
h_{t,x}^{(n)}\left( g\right) =h_{t,x}^{\left( \sqrt{c\log n}\right) }\left(
g\right) 
\]
For $h_{t,x}^{(n)}\in\mathcal{F}_{n}$ write 
\begin{equation*}
\alpha_{n}(h_{t,x}^{(n)})=\sum_{i=1}^{n}\frac{1\left\{ B_{i}(t)\leq
x,B_{i}\in\mathcal{C}_{n}\right\} -P\left\{ B(t)\leq x,B\in%
\mathcal{C}_{n}\right\} }{\sqrt{n}}.
\end{equation*}
Using (\ref{v}), note that for each $\left( t,x\right) \in\mathcal{T}\left(
\gamma_{n}\right) $, when $B_{i}\in\mathcal{C}_{n}$, for $i=1,\dots ,n$, 
\begin{equation*}
\begin{split}
\alpha_{n}\left( h_{t,x}^{(n)}\right) & =v_{n}\left( t,x\right) +\sqrt{n}%
P\left\{ B\left( t\right) \leq x,B\notin \mathcal{C}_{n}\right\}\\
& =\alpha_{n}\left( h_{t,x}\right) +\sqrt{n}P\left\{ B\left( t\right) \leq 
x,B\notin\mathcal{C}_{n}\right\} .  
\end{split}
\end{equation*}

Set 
\begin{equation*}
\mathcal{F}_{n}\left( \varepsilon\right) =\left\{ \left( f,f^{\prime
}\right) \in\mathcal{F}_{n}^{2}:d_{P}\left( f,f^{\prime}\right)
<\varepsilon\right\}
\end{equation*}
and 
\begin{equation*}
\mathcal{G}_{n}\left( \varepsilon\right) =\left\{ f-f^{\prime}:\left(
f,f^{\prime}\right) \in\mathcal{F}_{n}\left( \varepsilon\right) \right\} ,
\end{equation*}
where 
\begin{equation*}
d_{P}\left( f,f^{\prime}\right) =\sqrt{E\left( f\left( B\right)
-f^{\prime}\left( B\right) \right) ^{2}}.
\end{equation*}
By the arguments given in the Appendix of Kevei and Mason \cite{KM} the
classes $\mathcal{F}_{n}\left( \varepsilon\right) $ and $\mathcal{G}%
_{n}\left( \varepsilon\right) $ are \textit{pointwise measurable}. This means that the 
use of Talagrand's inequality below is justified.

Fix $n\geq1$. Let $B_1,\ldots,B_n$ be i.i.d.~$%
B$, and $\epsilon_{1},\ldots,\epsilon_{n}$ be independent Rademacher
random variables mutually independent of $B_{1},\ldots ,B_{n}$.
Write for $\varepsilon>0$, 
\begin{equation*}
\mu_{n}^{S}\left( \varepsilon\right) =E\left\{ \sup_{f-f^{\prime}\in\mathcal{%
G}_{n}\left( \varepsilon\right) }\left\vert \frac{1}{\sqrt{n}}%
\sum_{i=1}^{n}\epsilon_{i}\left( f-f^{\prime}\right) \left(
B_{i}\right) \right\vert \right\} ,
\end{equation*}

Observe that as long as $\varepsilon=\varepsilon_{n}$ and $\gamma=\gamma_{n}$
satisfy%
\begin{equation}
\sqrt{n}\varepsilon_{n}/\sqrt{\log n}\rightarrow\infty  \label{logn}
\end{equation}
and%
\begin{equation}
\log\left( \frac{\log n}{\varepsilon_{n}\gamma_{n}}\right) /\log
n\rightarrow\varsigma>0\text{, as }n\rightarrow\infty,  \label{eta2}
\end{equation}
we have 
\begin{equation*}
\sqrt{n}\varepsilon_{n}/\sqrt{\log\left( \frac{\log n}{\varepsilon_{n}%
\gamma_{n}}\right) }\rightarrow\infty,\text{ as }n\rightarrow\infty,
\end{equation*}
which by (57) in \cite{KM} implies
that for all large enough $n$ for a suitable $A_{1}>0$ 
\begin{equation*}
\mu_{n}^{S}\left( \varepsilon_{n}\right) \leq A_{1}\varepsilon _{n}%
\sqrt{\log\left( \frac{\log n}{\varepsilon_{n}\gamma_{n}}\right) }.
\end{equation*}
This, in turn, by (\ref{eta2}) gives for all large enough $n$, for some $A_1'>0$ 
\begin{equation*}
\mu_{n}^{S}\left( \varepsilon_{n}\right) \leq A_{1}^{\prime}\varepsilon_{n}\sqrt{%
\log n}.
\end{equation*}
Therefore by Talagrand's inequality (\ref{tal}) applied with $M=1$, we have
for suitable finite positive constants $D_1$, $D_1^{\prime}$, $D_20$ and for all $z>0$,
\begin{equation} \label{eq:ineq1}
\begin{split}
P\left\{ ||\sqrt{n}\alpha_{n}||_{\mathcal{G}_{n}\left( \varepsilon
_{n}\right) }\geq D_{1}^{\prime}(\varepsilon_{n}\sqrt{n\log n}+ z)\right\} &
\leq P\left\{ ||\sqrt{n}\alpha_{n}||_{\mathcal{G}_{n}\left(
\varepsilon_{n}\right) }\geq D_{1}(\sqrt{n}\mu_{n}^{S}\left( \varepsilon
_{n}\right) +z)\right\} \\
& \leq2\left\{ \exp\left( -\frac{D_{2} z^{2}}{n\sigma_{\mathcal{G}_{n}\left(
\varepsilon_{n}\right) }^{2}}\right) +\exp(-D_{2} z)\right\} .
\end{split}%
\end{equation}
Let 
\begin{equation}
\varepsilon_{n}=c_{1}\gamma_{n}^{-H/2}\left( \log\log n/n\right)
^{1/4},\text{ for some }c_{1}>0.  \label{epsi}
\end{equation}
Recall that $\gamma_{n}$ satisfies (\ref{eta}) with $\eta < 1/(2 H)$, which implies 
$\varepsilon_n \to 0$. Further, $\varepsilon_{n}$ fulfills (\ref{logn}) and 
\begin{equation*}
\frac{\log\left( \frac{\log n}{\varepsilon_{n}\gamma_{n}}\right) }{\log n}%
\rightarrow\frac{1}{4}+\eta\left( 1-\frac{H}{2}\right) =:\varsigma>0,
\end{equation*}
which says that (\ref{eta2}) holds. Also 
\begin{equation*}
n\sigma_{\mathcal{G}_{n}\left( \varepsilon_{n}\right) }^{2}=n\sup _{g\in%
\mathcal{G}_{n}\left( \varepsilon_{n}\right) }\mathrm{Var}(g(B))\leq
n\varepsilon_{n}^{2}.
\end{equation*}
Hence, 
\begin{equation*}
2\left\{ \exp\left( -\frac{D_{2} z^{2}}{n\sigma_{\mathcal{G}_{n}\left(
\varepsilon_{n}\right) }^{2}}\right) +\exp\left( -D_{2} z\right) \right\}
\leq 2\left\{ \exp\left( -\frac{D_{2} z^{2}}{n\varepsilon_{n}^{2}}\right)
+\exp\left( -D_{2} z\right) \right\} ,
\end{equation*}
which, with $z= \varepsilon_{n}\sqrt{dn\log n/D_{2}}$ for some $d>0$, is 
\begin{equation*}
\leq2\left\{ \exp\left( -d\log n\right) +\exp\left( -\sqrt{dD_{2}}%
\varepsilon_{n}\sqrt{n\log n}\right) \right\} .
\end{equation*}
By choosing $d>0$ large enough, (\ref{eq:ineq1}) combined with the
Borel--Cantelli lemma gives that, w.p.~1, 
\begin{align*}
||\alpha_{n}||_{\mathcal{G}_{n}\left( \varepsilon_{n}\right) } &
=\sup\left\{ \left\vert \alpha_{n}\left( h_{s,x}^{\left( n\right)
}-h_{t,y}^{\left( n\right) }\right) \right\vert :\left( s,x\right) ,\left(
t,y\right) \in\mathcal{T}\left( \gamma_{n}\right) \text{, }d_{P}^{2}\left(
h_{s,x}^{\left( n\right) },h_{t,y}^{\left( n\right) }\right)
<\varepsilon_{n}^{2}\right\} \\
& =O\left( n^{-1/4}\gamma_{n}^{-H/2}\left( \log\log n\right)
^{1/4}\left( \log n\right) ^{1/2}\right) .
\end{align*}
Recall that $\mathcal{T}\left( \gamma_{n}\right) =\left[ \gamma _{n},T%
\right] \times\mathbb{R}$. Since for $\gamma_{n}\leq t\leq T$%
\begin{equation*}
d_{P}^{2}\left( h_{t,x}^{\left( n\right) },h_{t,y}^{\left( n\right) }\right)
=E\left[ \left( 1\left\{ B(t) \leq x\right\}
-1\left\{ B(t) \leq y\right\} \right)
1\left\{ B\in\mathcal{C}_{n}\right\} \right] ^{2}
\end{equation*}%
\begin{equation*}
\leq E\left\vert 1\left\{ B(t) \leq x\right\}
-1\left\{ B(t) \leq y\right\} \right\vert
=\left\vert F\left( t,x\right) -F\left( t,y\right) \right\vert
\leq\gamma_{n}^{-H}\left\vert x-y\right\vert ,
\end{equation*}
i.e.~$|x-y| \leq c_1^2 \sqrt{(\log \log n) / n}$ implies that
$h_{t,x}^{(n)} - h_{t,y}^{(n)} \in \mathcal{G}_n(\varepsilon_n)$.
This says that, w.p.~$1$, with $c_{1}$ as in (\ref{epsi}),%
\begin{equation*}
\sup\left\{ \left\vert \alpha_{n}\left( h_{t,x}^{\left( n\right)
}-h_{t,y}^{\left( n\right) }\right) \right\vert :t\in\left[ \gamma _{n},T%
\right] ,\left\vert x-y\right\vert <\frac{c_{1}^{2}\sqrt{\log\log n}}{\sqrt{n%
}}\right\} \leq||\alpha_{n}||_{\mathcal{G}_{n}\left(
\varepsilon_{n}\right) ,}
\end{equation*}
where w.p.~$1$, 
\begin{equation*}
||\alpha_{n}||_{\mathcal{G}_{n}\left( \varepsilon_{n}\right) }=O\left(
n^{-1/4}\gamma_{n}^{-H/2}\left( \log\log n\right) ^{1/4}\left( \log
n\right) ^{1/2}\right) .
\end{equation*}
Next note that 
\begin{equation*}
\begin{split}
\Lambda_{n} & :=\sup\left\{ \left\vert \alpha_{n}\left( h_{t,x}\right)
-\alpha_{n}\left( h_{t,x}^{\left( n\right) }\right) \right\vert :\left(
t,x\right) \in\mathcal{T}\left( \gamma_{n}\right) \right\} \\
& \leq\sqrt{n}\sum_{i=1}^{n}1\left\{ B_{i}\notin\mathcal{C%
}_{n}\right\} +\sqrt{n}P\left\{ B\notin\mathcal{C}%
_{n}\right\} .
\end{split}%
\end{equation*}
We readily get using inequality (\ref{FLS1}) that for any $\omega>2$ there
exists a $c>0$ in (\ref{cn}) such that $P\left\{ B\notin \mathcal{C}_{n}\right\} \leq 
n^{-\omega},$ which implies 
\begin{equation*}
P\left\{ \Lambda_{n}>\sqrt{n}n^{-\omega}\right\} \leq n^{1-\omega}.
\end{equation*}
Thus we readily see by using the Borel--Cantelli lemma that, w.p.~$1$, 
\begin{equation} \label{an}
\begin{split}
& \sup\left\{ \left\vert \alpha_{n}\left( h_{t,x}-h_{t,y}\right) \right\vert
:t\in\left[ \gamma_{n},T\right] ,\left\vert x-y\right\vert <\frac{c_{1}^{2}%
\sqrt{\log\log n}}{\sqrt{n}} \right\}   \\
& =O\left( n^{-1/4}\gamma_{n}^{-H/2}\left( \log\log n\right)
^{1/4}\left( \log n\right) ^{1/2}\right) .  
\end{split}
\end{equation}
Applying Proposition \ref{prop:4} with $\delta= H/4$, keeping (\ref{less}) in 
mind, and by
choosing $c_{1}>0$ large enough in the definition of $\varepsilon_{n}$, we
see that, w.p.~$1,$ for all large $n$ 
\begin{equation*}
\sup_{\left( \alpha,t\right) \in\left[ \rho,1-\rho\right] \times\left[
\gamma_{n},T\right] }\left\vert \tau_{\alpha}^{n}\left( t\right)
-\tau_{\alpha}\left( t\right) \right\vert
\leq \frac{T^{3H/4} D_{0}\sqrt{\log\log n}}{%
\sqrt{n}}\leq\frac{c_{1}^{2}\sqrt{\log\log n}}{\sqrt{n}},
\end{equation*}
which says that, w.p.~$1,$ for all large enough $n$ uniformly in $\left(
\alpha,t\right) \in\left[ \rho,1-\rho\right] \times\left[ \gamma _{n},T%
\right] $, 
\[
\begin{split}
& \sup\left\{ \left\vert v_{n}\left( t,\tau_{\alpha}\left( t\right) \right)
-v_{n}\left( t,\tau_{\alpha}^{n}\left( t\right) \right) \right\vert :\left(
\alpha,t\right) \in\left[ \rho,1-\rho\right] \times\left[ \gamma_{n},T\right]
\right\} \\
& \leq\sup\left\{ \left\vert \alpha_{n}\left( h_{t,x}-h_{t,y}\right)
\right\vert :t\in\left[ \gamma_{n},T\right] ,\left\vert x-y\right\vert <%
\frac{c_{1}^{2}\sqrt{\log\log n}}{\sqrt{n}}\right\} .
\end{split}
\]
Thus by (\ref{an}), w.p.~$1$, for large enough $c>0$ and $c_{1}>0$, 
\begin{align*}
& \sup\left\{ \left\vert v_{n}\left( t,\tau_{\alpha}\left( t\right) \right)
-v_{n}\left( t,\tau_{\alpha}^{n}\left( t\right) \right) \right\vert :\left(
\alpha,t\right) \in\left[ \rho,1-\rho\right] \times\left[ \gamma_{n},T\right]
\right\} \\
& =O\left( n^{-1/4}\gamma_{n}^{-H/2}\left( \log\log n\right)
^{1/4}\left( \log n\right) ^{1/2}\right) .
\end{align*}
On account of (\ref{BHD}) this finishes the proof of Theorem \ref{th:3}.
\hfill$\square$

\subsection{Proof of Corollary \protect\ref{cor1}}

Let $\gamma _{n}=n^{-\eta}$, where $0<\eta <1/(2H)$ to be determined later. By Theorem 
\ref{th:3} 
\begin{equation} \label{C1}
\begin{split}
& \sup_{\left( t,\alpha \right) \in \left[ \gamma _{n},T\right] \times \left[
\rho ,1-\rho \right] }\left\vert t^{H}v_{n}\left( t,\tau _{\alpha }\left(
t\right) \right) +\frac{\exp \left( -\frac{z_{\alpha }^{2}}{2}\right) }{%
\sqrt{2\pi }}u_{n}\left( t,\alpha \right) \right\vert  \\
& =O\left( \left( \log \log n\right) ^{1/4}\left( \log n\right) ^{1/2}n^{-%
\frac{1}{4}+\eta  \frac{H}{2}  }\right) ,\text{a.s.}
\end{split}
\end{equation}
Next 
\begin{equation}
\sup_{\left( t,\alpha \right) \in \left[ 0,\gamma _{n}\right] \times \left[
\rho ,1-\rho \right] }\left\vert t^{H}v_{n}\left( t,\tau _{\alpha }\left(
t\right) \right) \right\vert \leq \sup \left\{ |t^{H}\alpha _{n}\left(
h_{t,x}\right) |:\left( t,x\right) \in \left[ 0,\gamma _{n}\right] \times 
\mathbb{R}\right\} .  \label{C11}
\end{equation}%
Now by a simple Borel--Cantelli argument based on inequality (\ref{max}) the
right side of (\ref{C11}) is equal to 
\begin{equation}
=O\left( \left( \log n\right) ^{1/2}n^{-\eta H}\right) ,\text{ a.s.}  \label{C2}
\end{equation}%
Next, by Proposition \ref{prop:4}, for any $0<\delta <H$ 
\begin{equation*}
\sup \left\{ \left\vert u_{n}\left( t,\alpha \right) \right\vert :\left(
\alpha ,t\right) \in \left[ \rho ,1-\rho \right] \times \left( a_{n},\gamma
_{n}\right] \right\} =O\left( \left( \log \log n\right) ^{1/2}n^{-\eta
\left( H-\delta \right) }\right) ,\text{ a.s.} 
\end{equation*}
so the same holds without the logarithmic factor
\begin{equation}
\sup \left\{ \left\vert u_{n}\left( t,\alpha \right) \right\vert :\left(
\alpha ,t\right) \in \left[ \rho ,1-\rho \right] \times \left( a_{n},\gamma
_{n}\right] \right\} =O\left( n^{-\eta \left( H-\delta \right) }\right) ,\text{ a.s.}  
\label{C3b}
\end{equation}
The latter rate is larger than the one in (\ref{C2}). Furthermore, comparing 
(\ref{C1}) and (\ref{C3b}) one sees that the optimal choice for $\eta$ is
$\eta =  1/(6H)$, and the best rate is $n^{-1/6 + \delta}$ (due to the arbitrariness of 
$\delta$ in the last step we changed $\eta \delta$ to $\delta$).
That is the statement is proved for $t\geq a_{n}$. 

Now we handle the $t < a_{n}$ case. Put 
\begin{equation*}
\Delta_{n}\left( a_{n}\right) :=\sup\left\{ \left\vert u_{n}\left(
t,\alpha\right) \right\vert :0\leq t\leq a_{n},0<\rho\leq\alpha\leq
1-\rho\right\} .
\end{equation*}
Observe that for all $0\leq t\leq a_{n}$ and $0<\rho\leq\alpha\leq1-\rho$, 
\begin{equation*}
\left\vert \tau_{\alpha}^{n}\left( t\right) \right\vert \leq\max_{1\leq
i\leq n}M_{i}\left( a_{n}\right) ,
\end{equation*}
where for $1\leq i\leq n$, $M_{i}\left( a_{n}\right) =\sup\left\{
\left\vert B_{i}\left( a_{n}s\right) \right\vert :0\leq s\leq1\right\} $. Notice that 
$B(a_{n}s) ,0\leq s\leq1,$ is equal in distribution to $a_{n}^{H} B( s)$, $0\leq s\leq1$. 
Further, as an application of the Landau--Shepp theorem we have for some $c>0$ and $d>0$%
\begin{equation}
P\left\{ \sup_{0\leq s\leq1}\left\vert B \left( s\right) \right\vert
>y\right\} \leq d\exp\left( -cy^{2}\right) ,\text{ for all }y>0.
\label{eq:max-est}
\end{equation}
We get now using a simple Borel--Cantelli lemma argument based on inequality
(\ref{eq:max-est}) and 
\begin{equation*}
M_{1}\left( a_{n}\right) \overset{\mathrm{D}}{=}a_{n}^{H}\sup_{{}}\left\{
\left\vert B_{1} ( s ) \right\vert :0\leq s\leq1\right\} ,
\end{equation*}
that for some $D>0$, w.p.~$1$, for all $n$ sufficiently large, 
\begin{equation*}
\max_{1\leq i\leq n}M_{i}\left( a_{n}\right) \leq Da_{n}^{H}\sqrt{\log n}.
\end{equation*}
Hence, w.p.~1, uniformly $0\leq t\leq a_{n}$ and $0<\rho\leq\alpha\leq1-\rho$, for all 
large enough $n$, 
\begin{equation*}
\left\vert \tau_{\alpha}^{n}\left( t\right) \right\vert \leq Da_{n}^{H}\sqrt{%
\log n}.
\end{equation*}
Also trivially we have uniformly $0\leq t\leq a_{n}$ and $0<\rho\leq\alpha
\leq1-\rho$ 
\begin{equation*}
\left\vert \tau_{\alpha}\left( t\right) \right\vert =t^{H}|z_{\alpha}|\leq
a_{n}^{H}z_{1-\rho}.
\end{equation*}
Thus, w.p.~1, uniformly $0\leq t\leq a_{n}$ and $0<\rho\leq\alpha\leq1-\rho$%
, for all large enough $n$, 
\begin{equation*}
\Delta_{n}\left( a_{n}\right) \leq2Da_{n}^{H}\sqrt{n\log n}.
\end{equation*}
Note that 
\begin{equation*}
\rho_{n}:=2Da_{n}^{H}\sqrt{n\log n}=2DC^{H}\left( \frac{\log\log n}{n}%
\right) ^{H/(2\delta)}\sqrt{n\log n}\text{,}
\end{equation*}
satisfies 
\begin{equation*}
\frac{-\log\rho_{n}}{\log n}\rightarrow\frac{H}{2\delta}-\frac{1}{2}=\frac{%
H-\delta}{2\delta}>0.
\end{equation*}
For some $\delta > 0$ small enough $(H - \delta)/ (2 \delta) > 1/6$, therefore
\[
\Delta_n(a_n) = O( n^{-1/6}), \text{ a.s.} 
\]
which together with (\ref{C1}), and (\ref{C2}) finish the proof
of the corollary. \hfill$\square$ \smallskip

\section{Appendix: Useful inequalities}

\subsection{Talagrand's inequality}

We shall be using the following exponential inequality due to Talagrand \cite%
{Talagrand}. \medskip

\noindent\textbf{Talagrand Inequality.} \textit{Let $\mathcal{G}$ be a
pointwise measurable class of measurable real-valued functions defined on a
measure space $(S,\mathcal{S})$ satisfying $||g||_{\infty}\leq M,\ g\in 
\mathcal{G}$, for some $0<M<\infty$. Let $X,X_{n}$, $n\geq1$, be a sequence
of i.i.d.~random variables defined on a probability space $\left( \Omega ,%
\mathcal{A},P\right) $ and taking values in $S$, then for all $z>0$ we have
for suitable finite constants $D_{1},D_{2}>0$, 
\begin{equation}
\begin{split}
P\left\{ ||\sqrt{n}\alpha_{n}||_{\mathcal{G}}\geq
D_{1}\left( E\left\Vert \sum_{i=1}^{n}\epsilon_{i}g(X_{i})\right\Vert _{%
\mathcal{G}}+ z\right) \right\} 
\leq 2\exp\left( -\frac{D_{2} z^{2}}{n\sigma_{\mathcal{G}}^{2}}\right)
+2\exp\left( -\frac{D_{2}z}{M}\right) ,
\end{split}
\label{tal}
\end{equation}
where $\sigma_{\mathcal{G}}^{2}=\sup_{g\in\mathcal{G}}\mathrm{Var}(g(X))$
and $\epsilon_{n},$ $n\geq1$, are independent Rademacher random
variables mutually independent of $X_{n}$, $n\geq1$.
}

\subsection{Application of Landau--Shepp Theorem}

By the L\'{e}vy modulus of continuity theorem for fractional Brownian motion 
$B^{(H)}$ with Hurst index $0<H<1$ (see Corollary 1.1 of Wang \cite{Wang}),
we have for any $0<T<\infty$, w.p.~$1$,%
\begin{equation*}
\sup_{0\leq s\leq t\leq T}\frac{\left\vert B^{(H)}\left( t\right)
-B^{(H)}\left( s\right) \right\vert }{f_{H}\left( t-s\right) }=:L<\infty.
\end{equation*}
Therefore we can apply the Landau and Shepp \cite{LandauShepp} theorem (also
see Sato \cite{Sato} and Proposition A.2.3 in \cite{VaartWellner}) to infer
that for appropriate constants $C>0$ and $D>0$, for all $z>0$,
\begin{equation}
P\left\{ L> z\right\} \leq C\exp\left( -D z^{2}\right) .  \label{FLS1}
\end{equation}

\subsection{A maximal inequality}

The following inequality is proved in Kevei and Mason \cite{KM}, where it is
Inequality 2.\medskip

\noindent \textbf{Inequality}
\textit{For all  $0<\gamma \leq 1$ and $\tau >0$ we have for some $E(\tau )$ and for 
suitable finite positive constants $D_{3},D_{4}>0$, for all $z>0$ 
\begin{equation}
\begin{split}
& P\left\{ \max_{1 \leq m \leq n} \sup_{(t,x)\in [ 0,\gamma ]\times 
\mathbb{R}}|\sqrt{m} t^{\tau } \alpha_{m} \left( h_{t,x}\right) |\geq
D_{3} \sqrt{n} \left( E(\tau ) (2 \gamma)^{\tau } +z\right) \right\} \\
& \leq 2\left\{ \exp \left( - D_{4} z^{2} \left( 2\gamma \right)^{-2\tau }\right) +\exp 
\left( -D_4 z \sqrt{n} \left( 2\gamma \right) ^{-\tau
}\right) \right\} .
\end{split}
\label{max}
\end{equation}
}

%
%

\end{document}